\documentclass[a4paper,12pt,twoside]{article}
\usepackage{amsbsy}
\usepackage{amssymb}
\usepackage{amsfonts}
\usepackage{amsmath}
\usepackage{amsopn}
\usepackage{amsthm}
\usepackage[english]{babel}
\usepackage[T1]{fontenc}

\newtheorem{teo}{Theorem}
\newtheorem{lemma}{Lemma}
\newtheorem{prop}{Proposition}
\newtheorem{cor}{Corollary}

\theoremstyle{definition}
\newtheorem{remark}{Remark}

\theoremstyle{definition}
\newtheorem{example}{Example}
\newtheorem*{gracias}{Acknowledgements}

\theoremstyle{remark}

\DeclareMathOperator{\er}{\mathbf{E}}
\DeclareMathOperator{\pr}{\mathbf{P}}

\DeclareMathOperator{\re}{\mathbb{R}}

\DeclareMathOperator{\van}{\xrightarrow[n\to\infty]{}}
\DeclareMathOperator{\vat}{\xrightarrow[t\to\infty]{}}

\title{Sina\v\i 's condition for real valued L\'evy processes}

\author{\textbf{V\'\i ctor RIVERO} \thanks{\'Equipe MODAL'X, Universit\'e Paris X Nanterre, UFR SEGMI, 200 Avuenue de la R\'epublique, 92000 Nanterre CEDEX France; E--mail: rivero@ccr.jussieu.fr}}
\date{\today}
\begin{document}

\maketitle

\begin{abstract}
We prove that the upward ladder height subordinator $H$ associated
to a real valued L\'evy process $\xi$ has Laplace exponent
$\varphi$ that varies regularly at $\infty$ (resp. at $0$) if and
only if the underlying L\'evy process $\xi$ satisfies Sina\v\i 's
condition at $0$ (resp. at $\infty$). Sina\v\i 's condition for
real valued L\'evy processes is the continuous time analogue of
Sina\v\i 's condition for random walks. We provide several
criteria in terms of the characteristics of $\xi$ to determine
whether or not it satisfies Sina\v\i 's condition. Some of these
criteria are deduced from tail estimates of the L\'evy measure of
$H,$ here obtained, and which are analogous to the estimates of
the tail distribution of the ladder height random variable of a
random walk which are due to Veraverbeke and Gr\"ubel.
\end{abstract}

\noindent\textbf{Key words:} L\'evy processes, Fluctuation theory, Regular Variation, long tailed L\'evy measures.

\noindent\textbf{MSC:} 60 G 30 (60 G 51).

\begin{section}{Introduction and main result}

Let $\xi=\{\xi_t,t\geq 0\}$ be a real valued L\'evy process,
$S=(S_t,t\geq 0)$ its current supremum and $L=(L_t, t\geq 0)$ the
local time at 0 of the strong Markov process $\xi$ reflected at
its current supremum, that is to say $(S_t-\xi_t,t\geq 0).$ In
this work we will obtain some asymptotic properties of the
ascending ladder height subordinator $H$ associated to $\xi$ (that
is, the current supremum of $\xi$ evaluated at the inverse of the
local time at 0, i.e. $L^{-1},$ $H\equiv(S_{L^{-1}_t},t\geq 0)$).
According to Fristedt~\cite{Fristedt74} the ascending ladder
process $(L^{-1},H)$ is a bivariate subordinator, that is, a
L\'evy process in $\re^2$ with increasing paths (coordinatewise)
whose bivariate Laplace exponent $\kappa,$
$$e^{-\kappa(\lambda_1,\lambda_2)}\equiv \er(e^{-\lambda_1L^{-1}_1-\lambda_2H_1}),\qquad \lambda_1, \lambda_2 \geq 0,$$ with the assumption $e^{-\infty}=0,$ is given by
$$\kappa(\lambda_1,\lambda_2)=\mathrm{k}\exp\left\{\int^{\infty}_0\frac{\mathrm{d}t}{t}\int_{[0,\infty[}(e^{-t}-e^{-\lambda_1t - \lambda_2x})\pr(\xi_t\in \mathrm{d}x)\right\},\qquad \lambda_1, \lambda_2 \geq 0,$$
with $\mathrm{k}$ a constant that depends on the normalization of
the local time. (See Doney~\cite{doney2001}, for a survey, and
Bertoin~\cite{MR98e:60117}~VI, for a detailed exposition of the
fluctuation theory of L\'evy processes and
Vigon~\cite{MR2002i:60101} for a description of the L\'evy measure
of $H$.)

The fact that the ladder process $(L^{-1},H)$ is a bivariate
subordinator is central in the fluctuation theory of L\'evy
processes because it enables to obtain several properties of the
underlying L\'evy process using results for subordinators, which
are objects simpler to manipulate. Among the various properties
that can be obtained using this fact, there is a well known
arc-sine law in the time scale for L\'evy processes, see Theorem
VI.3.14 in Bertoin's book~\cite{MR98e:60117} for a precise
statement. That result tell us that Spitzer's condition is a
condition about the underlying L\'evy process $\xi$ which ensures
that the Laplace exponent $\kappa(\cdot,0)$ of the ladder time
subordinator $L^{-1}$ is regularly varying and which in turn
permits to obtain an arc-sine law in the time scale for L\'evy
processes. Now, if we want to establish an analogous result in the
space scale we have to answer the question: \textit{What is the
analogue of Spitzer's condition for the upward ladder height
process $H$?} or put another way: \textit{What do we need to
assume about $\xi$ to ensure that the Laplace exponent
$$\varphi(\lambda)\equiv\kappa(0,\lambda)=\mathrm{k}\exp\left\{\int^{\infty}_0\frac{\mathrm{d}t}{t}\int_{[0,\infty[}(e^{-t}-e^{-\lambda x})\pr(\xi_t\in \mathrm{d}x)\right\},\qquad \lambda \geq 0,$$
of $H$ varies regularly}?

To motivate the answer that we will provide to these questions we
will make a slight digression to recall the analogue of this
problem for random walks. Let $X_1, X_2,\ldots $ be a sequence of
independent and identically distributed random variables and $Z$
its associated random walk $Z_0=0,\ Z_n=\sum^n_{k=1} X_k, n>0.$
The analogue of the upward ladder process for $Z$ is given by $(N,
Z_N)$ where $N$ is the first ladder epoch of the random walk $Z,$
$N=\min\{k: Z_k>0\},$ and $Z_N,$ is the position of $Z$ at the
instant $N.$ We denote by $M_n=\sup\{Z_k, 0\leq k\leq n\}, n\geq
0$ the current supremum of the random walk and by $T(x)$ the first
instant at which the random walk passes above the level $x$,
$T(x)=\inf\{n>0: Z_n\geq x\}$.

Greennwood, Omey and Teugels~\cite{MR85e:60093} proved that the condition,
$$\qquad  \sum^\infty_{n=1}\frac{1}{n}\pr(z < Z_n \leq \lambda z)\xrightarrow[z\to\infty]{} \beta \log(\lambda),\qquad \forall \lambda >1,$$ for some $\beta\in[0,1],$
which they called Sina\v\i 's condition making reference to a work
of Sina\v\i 's~\cite{sinai57}, plays the same r\^ole, in the
obtaining of limit results for the ladder height, as does
Spitzer's condition for the ladder time. More precisely, when the
random walk does not drift to $-\infty,$ we have the following
equivalences.

\begin{teo}[Dynkin~\cite{dynkin55}, Rogozin~\cite{rogozin71}, Greenwood et al.~\cite{MR85e:60093}] Assume $\sum_n \frac{1}{n}\pr(Z_n>0)=\infty.$ For $\beta\in [0,1],$ the following are equivalent
\begin{enumerate}
\item[(i)] ${\displaystyle
\lim_{z\to\infty}\sum_{n}\frac{1}{n}\pr(z< Z_n \leq \lambda
z)=\beta \log(\lambda)}, \qquad \forall \lambda>1;$ \item[(ii)]
$\int^x_0\pr(Z_N > u)\mathrm{d}u \sim x^{1-\beta}l(x)$ as
$x\to\infty,$ with $l$ slowly varying; \item[(iii)] The random
variables $\gamma_x=\frac{\displaystyle M_{T(x)-1}}{\displaystyle
x}$ converge in distribution as $x \to \infty;$ the limit law is
the generalized arc-sine law of parameter $\beta,$ that is to say
that it  has a density
$$q_{\beta}(y)=\frac{\sin(\beta \pi)}{\pi}y^{-\beta}(1-y)^{\beta-1}, \qquad y\in ]0,1[,$$ if $\beta\in ]0,1[;$ and is degenerate with unit mass at 1 or
0, according as $\beta=1$ or $\beta=0.$
\end{enumerate}
\end{teo}
The equivalence between (i) and (ii) is due to Greenwood et
al.~\cite{MR85e:60093} (see also \cite{MR90i:26003} Theorem
8.9.17), that between (ii) and (iii) is due to
Dynkin~\cite{dynkin55} in the case $\beta\in]0,1[$ and to
Rogozin~\cite{rogozin71} in the other cases.

Thus, the previous theorem tells us that Sina\v\i 's condition
enables to obtain a spatial arc-sine law for random walks and,
given that the fluctuation theory for Levy processes mirrors that
of random walks, it is natural to hope that the answer to the
questions posed above is the continuous time version of Sina\v\i
's condition.

We will say that a L\'evy process $\xi$ satisfies Sina\v\i 's
condition at $\infty$ (resp. at $0$) if \begin{description}
\item[(Sina\v\i)] There exists a $0\leq \beta\leq 1$ such that
$$\int^{\infty}_{0}\frac{\mathrm{d}t}{t}\pr(z<\xi_t \leq\lambda
z)\longrightarrow\beta\log(\lambda) \quad \text{as} \quad z\to
+\infty \quad \text{(resp.  $z\to 0+$)}, \qquad \forall
\lambda>1.$$
\end{description} The term $\beta$ will be called Sina\v\i 's
index of $\xi.$ Observe that if Sina\v\i 's condition hold we have
that
$$\int^{\infty}_{0}\frac{\mathrm{d}t}{t}\pr(\lambda z<\xi_t \leq
z)\longrightarrow -\beta\log(\lambda) \quad \text{as} \quad z\to
+\infty \quad \text{(resp.  $z\to 0+$)},\qquad \forall
\lambda\in]0,1[.$$

\begin{example}
A L\'evy process, $\xi,$ which satisfies Sina\v\i 's condition is
the strictly stable process with index $0<\alpha\leq 2.$ Indeed,
for every $z>0$ and $\lambda>1$ we have by the scaling property of
$\xi$ that
\begin{equation*}
\begin{split}
\int^\infty_0\frac{\mathrm{d}t}{t}\pr\left(z<\xi_t\leq \lambda z\right)&=\int^\infty_0\frac{\mathrm{d}t}{t}\pr\left(z<t^{1/\alpha}\xi_1\leq\lambda z\right)\\
&=\er\left(1_{\{\xi_1>0\}}\int^\infty_0\frac{\mathrm{d}t}{t}1_{\{(z/\xi_1)^\alpha< t \leq (z\lambda/\xi_1)^\alpha\}}\right)\\
&=\er\left(1_{\{\xi_1>0\}}\log(\lambda^\alpha)\right)=\alpha\pr(\xi_1>0)\log(\lambda).\end{split}
\end{equation*}
Thus any stable process $\xi$ does satisfies Sina\v\i 's condition at infinity
and at 0 with index $\alpha\rho,$ where $\rho$ is the positivity
parameter of $\xi,$ $\rho=\pr(\xi_1\geq 0).$
\end{example}
We recall that a measurable function $f:[0,\infty[\to [0,\infty[$
varies regularly at infinity (resp. at 0) with index
$\alpha\in\re,$ $f\in RV^{\infty}_{\alpha}$ (resp. $\in
RV^{0}_{\alpha}$), if for any $\lambda>0,$ $$\lim \frac{f(\lambda
x)}{f(x)}=\lambda^{\alpha}\qquad \text{at}\ \infty\ \text{(resp.
at $0$)}.$$

We have all the elements to state our main result, which provides
an answer to the questions above.

\begin{teo}\label{mainthm:1}For $\beta\in[0,1],$ the following are equivalent
\begin{enumerate}
\item[(i)] The L\'evy process $\xi$ satisfies Sina\v\i 's  condition at $\infty$ (resp. at 0) with index $\beta.$
\item[(ii)] The Laplace exponent of the ladder height subordinator $H$ varies regularly at $0$ (resp. at $\infty$) with index $\beta.$
\end{enumerate}
\end{teo}

The proof of this result is an easy consequence of a fluctuation identity due to Bertoin and Doney~\cite{MR96b:60190} and Theorem~\ref{teo1:s} below.

\begin{proof}
By the fluctuation identity of Bertoin and Doney we have that for any $z>0,\lambda>1,$
$$\int^\infty_0\frac{\mathrm{d}t}{t}\pr(z< \xi_t \leq\lambda z)=\int^\infty_0\frac{\mathrm{d}t}{t}\pr(z< H_t \leq\lambda z).$$ As a consequence, Sina\v\i 's condition is satisfied by the L\'evy process $\xi$ if and only if it is satisfied by the ascending ladder height subordinator $H.$ The result then follows from Theorem~\ref{teo1:s}, which establishes that the Laplace exponent $\phi$ of any given subordinator, say $\sigma,$ varies regularly if and only if $\sigma$ satisfies Sina\v\i 's condition. \end{proof}

Assuming that the L\'evy process $\xi$ satisfies Sina\v\i 's
condition and applying known results for subordinators, when its
Laplace exponent is regularly varying, we can deduce the behavior
at 0 or $\infty$ of $\xi$ from that of $H.$  (See
Bertoin~\cite{MR98e:60117} Ch. III for an account on the short and
long time behavior of subordinators.) The following
\textit{spatial arc-sine law} for L\'evy processes is an example
of the results that can be obtained.

\begin{cor}\label{corsec1:arcsin} For $r>0,$ denote the first exit time of $\xi$ out of $]-\infty,r]$ by $T_r=\inf\{t>0: \xi_t>r\},$ the undershot and overshot of the supremum of $\xi$ by $U(r)=r-S_{T_r-}$ and $O(r)=S_{T_r}-r=\xi_{T_r}-r.$ For any $\beta\in [0,1],$ the conditions (i) and (ii) in Theorem~\ref{mainthm:1} are equivalent to the following conditions:
\begin{enumerate}
\item[(iii)] The random variables $r^{-1}(U(r), O(r))$ converge in
distribution as $r\to\infty$ (respectively, as $r\to 0$).
\item[(iv)] The random variables $r^{-1}O(r)$ converge in
distribution as $r\to\infty$ (respectively, as $r\to 0$).
\item[(v)] The random variables $r^{-1}S_{T_r-}$ converge in
distribution as $r\to\infty$ (respectively, as $r\to 0$).
\item[(vi)] $\lim r^{-1}\er\left(S_{T_r-}\right)=\beta\in[0,1]$ as
$r\to\infty$ (respectively, as $r\to 0$).
\end{enumerate}
In this case, the limit distribution in (iii) is determined as
follows: if $\beta=0$ (resp. $\beta=1$), it is the Dirac mass at
$(1,\infty)$ (resp. at $(0,0))$. For  $\beta\in]0,1[,$ it is the
distribution with density
$$p_{\beta}(u,w)=\frac{\beta\sin\beta\pi}{\pi}(1-u)^{\beta-1}(u+w)^{-1-\beta},\qquad
0<u<1, w>0.$$ In particular, the limit law in (v) is the
generalized arc-sine law of parameter $\beta.$
\end{cor}

\begin{proof}
We recall that for every $r>0,$ the random variables $(U(r),
O(r))$ are almost surely equal to the undershot and overshoot,
$(U_H(r), O_H(r)),$ of the ladder height subordinator $H.$ Thus
the result is a straightforward consequence of the Dynkin-Lamperti
arc-sine law for subordinators, Theorem~III.3.6
in~\cite{MR98e:60117}, using the elementary relations: for every
$r>0$
\begin{gather*}
\pr(U_H(r)>y)=\pr(O_H(r-y)>y),\qquad r>y>0,\\
\pr(O_H(r)>x, U_H(r)>y)=\pr(O_H(r-y)>x+y),\qquad r>y>0,\ x>0.
\end{gather*}
\end{proof}

We would like to remark that in the case where Sina\v\i 's
condition hold for $\xi$ with index $\beta=1$ at $\infty$ (resp.
at 0) we have that $r^{-1}\xi_{T_r}$ converges in law to $1$ as
$r\to \infty$ (resp. to $0$). The almost sure convergence of this
random variable was studied by Doney and
Maller~\cite{MR2003d:60090}. Precisely, Theorem 8 of Doney and
Maller~\cite{MR2003d:60090}  provide necessary and sufficient
conditions, on the characteristics of $\xi,$ according to which
$r^{-1}\xi_{T_r}$ converges a.s. to 1 as $r\to\infty.$ Moreover,
Theorem 4 in~\cite{MR2003d:60090} establish that the latter r.v.
converge a.s. to 1 as $r\to 0$ if and only if $\xi$ creeps upward.

To summarize, in Theorem~\ref{mainthm:1} we provided a necessary
and sufficient condition in terms of the marginal laws of $\xi$
which completely answers the questions posed at the beginning of
this work. However, the possible drawback of this result is that
in most of the cases we only know the characteristics of the
L\'evy process $\xi,$ that is, its linear and Gaussian terms and
L\'evy measure, and so it would be suitable to have a condition in
terms of the characteristics of the process. That is the purpose
of the second part of this work.

One case at which Sina\v\i 's condition can be verified using the
characteristics of the process is the case at which the underlying
L\'evy process belongs to the domain of attraction at infinity
(respectively, at 0) of a strictly stable law of index
$0<\alpha\leq 2,$ and which does not require a centering function.
That is, whenever there exists a deterministic function
$b:]0,\infty[\to]0,\infty[$ such that
\begin{equation}\label{UDS} \frac{\xi_t}{b(t)} \xrightarrow[]{\mathcal{D}}X(1), \quad\text{as}\ t\to\infty\ (\text{respectively, as $t\to 0$}),
\end{equation} with $X(1)$ a strictly stable random variable of parameter $0<\alpha\leq 2.$ It is well known that if such a function $b$ exists, it is regularly varying at infinity (respectively, at $0$) with index $\beta=1/\alpha$.
Plainly, the convergence in (\ref{UDS}) can be determined in terms
of the characteristic exponent $\Psi$ of $\xi,$ i.e.
$\er(e^{i\lambda\xi_t})=\exp\{t\Psi(\lambda)\},\lambda\in\re,$
since the latter convergence in distribution is equivalent to the
validity of the limit
\begin{equation}\label{eq3ch1:DA}\lim  t \Psi\left(\frac{\lambda}{b(t)}\right)=\Psi_{\alpha}(\lambda),\quad \text{as}\ t\to\infty,\quad \text{(respectively, as $t\to 0$ )}\ \text{for}\ \lambda\in\re,\end{equation}
where $\Psi_{\alpha}$ is the characteristic exponent of a strictly
stable law and is given by
\begin{equation*}
\Psi_{\alpha}(\lambda)=\begin{cases}-c|\lambda|^\alpha\left(1-i\delta\textup{sgn}(\lambda)\tan(\pi\alpha/2)\right)
& \ 0<\alpha<1\ \text{or}\ 1<\alpha<2;\\
-c|\lambda|^\alpha\left(1-i\delta\textup{sgn}(\lambda)\tan(\pi\alpha/2)\ln(|\lambda|)\right),
& \ \alpha=1;\\
-q^2\lambda^2/2, & \ \alpha=2;\end{cases}
\end{equation*}
for $\lambda\in\re,$ where $c>0$ and the term $\delta\in[-1,1]$ is
the so called skewness parameter. We have the following theorem
whose proof will be given in Section~\ref{proofthm2}.

\begin{teo}\label{teoch1:DA} Let $0<\alpha\leq 2$ and $\delta\in[-1,1].$ Assume that there exits a function $b:]0,\infty[\to]0,\infty[$ such that the limit in equation~(\ref{eq3ch1:DA}) holds as $t$ goes to infinity (resp. as $t\to 0$). Then the L\'evy process $\xi$ satisfies Sinai's condition at $\infty$ (respectively, at 0) with index $\alpha\rho,$ where $\rho$ is given by $\rho=1/2+(\pi\alpha)^{-1}\arctan(\delta\tan\left(\alpha\pi/2\right)).$
\end{teo}

With the aim of providing some other criteria in terms of the
characteristics of the underlying L\'evy process $\xi$ to
determine whether or not it satisfies Sina\v\i 's condition, in
Section~\ref{somecriterions} we will provide, under some assumptions, some estimates
of the L\'evy measure of $H.$ Those estimates are more general
than needed and are the L\'evy processes version of the results of
Veraverbeke~\cite{MR0423543} and Gr\"ubel~\cite{Grubel85} for
random walks. Furthermore, they are in the same vein as those
results obtained recently by Kl\"uppelberg, Kyprianou and
Maller~\cite{KKM}, Propositions 5.3 and 5.4.

Using those  estimates we will provide some criteria to determine
whether or not the L\'evy measure of $H$ is regularly varying. For
our ends this will be sufficient because by an Abelian-Tauberian
theorem we know that this determines the regular variation of the
Laplace exponent of $H.$

The rest of this note is organized as follows. In Section~2, we
state and prove an equivalent form of the Dynkin \& Lamperti's
theorem for subordinators (see e.g.~\cite{MR98e:60117}
Theorem~III.6) which is interesting in itself. In
Section~\ref{proofthm2}, we prove Theorem~\ref{teoch1:DA}.
Finally, in Section~\ref{somecriterions} we provide some estimates
of the L\'evy measure of $H$ and we use those estimates to provide
some criterions in terms of the characteristics of the underlying
L\'evy process $\xi$ to determine whether or not it satisfies
Sina\v\i 's condition. Section~\ref{somecriterions} is self
contained and can be read separately from the rest of the paper.
\end{section}

\begin{section}{A result for subordinators}
Let $\sigma=(\sigma_t,t\geq 0)$ be a subordinator, possibly
killed, with life time $\zeta,$ and denote by $\phi$ its Laplace
exponent, $$\phi(\lambda)\equiv -\log\er(e^{-\lambda \sigma_1},
1<\zeta),\qquad \lambda \geq 0.$$ It is well known that the
Laplace exponent $\phi$ can be represented as
$$\phi(\lambda)=\kappa + \lambda d +
\int_{]0,\infty[}(1-e^{-\lambda y})\nu(\mathrm{d}y),\qquad \lambda
\geq 0,$$ where $\kappa,d\geq 0$ are the killing rate and drift
coefficient of $\sigma,$ respectively, and $\nu$ is the L\'evy
measure of $\sigma$, that is, a measure on $]0,\infty[$ such that
$\int_{]0,\infty[} \min\{1,y\} \nu(\mathrm{d}y) <\infty.$

The aim of this section is prove the following equivalent form of
the Dynkin \& Lamperti's theorem for subordinators, see
e.g.~\cite{MR98e:60117} Theorem~III.6.

\begin{teo}\label{teo1:s}
For $\beta\in[0,1],$ the following are equivalent:
\begin{enumerate}
\item[(i)] The subordinator $\sigma$ satisfies Sina\v\i 's condition at $0+$ (resp. at $+\infty$) with index $\beta.$
\item[(ii)] The Laplace exponent $\phi$ is regularly varying at $+\infty$ (resp. at $0+$) with index $\beta.$
\end{enumerate}
\end{teo}

An account on necessary and sufficient conditions according to
which a subordinator has a Laplace exponent such that $\phi\in
RV^{\infty}_{\beta}$ or $\phi\in RV^{0}_{\beta}$ can be found in
\cite{MR98e:60117} page 82.

The  proof of this result relies on the following elementary remark.
\begin{remark}Write $$\phi(\theta)=(1+\phi(\theta))\frac{\phi(\theta)}{1+\phi(\theta)},\qquad \theta\geq 0.$$ The
first (resp. second) factor in the right hand term of the previous
equality can be used to determine the behavior at infinity (resp.
at 0) of $\phi.$ More precisely, $\phi\in RV^{\infty}_{\beta},$
(resp. $\in RV^{0}_{\beta}$) if and only  $$1+\phi(\cdot)\in
RV^{\infty}_{\beta}, \quad \text{(resp.}\
\frac{\phi(\cdot)}{1+\phi(\cdot)}\in RV^{0}_{\beta}\text{)}.$$
\end{remark}

The proof of Theorem~\ref{teo1:s} will be given using the previous
remark and via three lemmas whose proof will be given at the
end of this section. The first of them will enable us to relate
the factors of the latter remark with a transformation of the
type Mellin's convolution.

\begin{lemma}\label{lemmaML:s} We have that
\begin{enumerate}
\item[(i)] ${\displaystyle 1+\phi(\theta)= \exp\{
\widehat{G}_{1}(\theta)\}},$ for $\theta>0;$ where the function
$\widehat{G}_{1}$ is the Mellin convolution of the non-decreasing
function
$$G_{1}(y)=\int^{\infty}_0\frac{\mathrm{d}t}{t}e^{-t}\pr(\sigma_t>1/y),\qquad
y>0,$$ and the kernel $k(x)=x e^{-x}, x>0;$ that is,
$$\widehat{G}_{1}(\theta)=k^{M}*G_{1}(\theta):=\int^{\infty}_0\frac{\mathrm{d}x}{x}k(\theta/x)G_{1}(x), \qquad \theta>0.$$

\item[(ii)]${\displaystyle\frac{\phi(\theta)}{1+\phi(\theta)}=\exp\left\{-\widehat{G}_{2}(\theta)\right\}},$
 for $\theta>0;$ where $\widehat{G}_{2}$ is the Laplace transform of
the measure
$$G_{2}(\mathrm{d}x)=\int^{\infty}_0\frac{\mathrm{d}t}{t}(1-e^{-
t})\pr(\sigma_t\in \mathrm{d}x),\qquad x>0;$$ which is in fact the
harmonic renewal measure associated to the  law
$F(\mathrm{d}x)=\pr(\sigma_{\Theta}\in \mathrm{d}x)$ with $\Theta$
an independent random variable with exponential law of parameter
$1$.
\end{enumerate}
\end{lemma}

A consequence of Lemma~\ref{lemmaML:s} is that  $1+\phi\in
RV^{\infty}_{\beta}$ if and only if
\begin{equation}\label{dif1:s}\lim_{\theta\to\infty}\widehat{G}_{1}(\lambda
\theta)-\widehat{G}_{1}(\theta)=\beta\log\lambda, \qquad
\forall\lambda>0.\end{equation} Moreover,
$\phi(\cdot)/(1+\phi(\cdot) )\in RV^{0}_{\beta}$ if and only if
\begin{equation}\label{dif2:s}\lim_{\theta\to
0}\widehat{G}_{2}(\lambda
\theta)-\widehat{G}_{2}(\theta)=-\beta\log\lambda, \qquad
\forall\lambda>0.\end{equation}

The second of these Lemmas enable us to relate Sinai's condition
with the behavior at infinity of the differences of the function
$G_{1},$ and those of the function $G_2(x)\equiv G_2[0,x],$ $x>0.$

\begin{lemma}\label{lemmast:1} Let $\beta\in[0,1].$
\begin{enumerate}
\item[(i)] Sina\v\i 's condition holds at 0 with index $\beta$ if and only if
\begin{equation*}%\label{eq2:s}
\lim_{z\to \infty}G_{1}(\lambda
z)-G_{1}(z)=\beta\log(\lambda),\qquad \forall
\lambda>1.\end{equation*}

\item[(ii)] Let $G_{2}(z):=G_{2}[0,z], z>0.$ Sina\v\i 's condition
holds at infinity with index $\beta$ if and only if
\begin{equation*}%\label{eq2:s}
\lim_{z\to \infty}G_{2}(\lambda
z)-G_{2}(z)=\beta\log(\lambda),\qquad \forall
\lambda>1.\end{equation*}
\end{enumerate}
\end{lemma}

The last ingredient to achieve the proof of Theorem~\ref{teo1:s}
is an Abelian--Tauberian's type result relying the behavior of the
differences of $G_{1}$ (resp. $G_{2}$) with those of the functions
$\widehat{G}_{1}$ (resp. $\widehat{G}_{2}$).

\begin{lemma}\label{lemmaAT:1}
\begin{enumerate}\item[(i)]The following are equivalent
\begin{equation}\label{eq2:s}
\begin{split}
\lim_{y\to \infty}G_{1}(\lambda
y)-G_{1}(y)&=\beta\log(\lambda),\qquad \forall
\lambda>0.\end{split}
\end{equation}

\begin{equation}\label{eq2e:s}
\lim_{\theta\to\infty}\widehat{G}_{1}(\lambda
\theta)-\widehat{G}_{1}(\theta)= \beta\log(\lambda), \qquad
\forall \lambda>0.
\end{equation}
Both imply that
$$G_{1}(\theta)-\widehat{G}_{1}(\theta)\xrightarrow[\theta\to\infty]{}\beta \gamma.$$

\item[(ii)]The following are equivalent
\begin{equation}\label{eq22:s}
\begin{split}
\lim_{y\to \infty}G_{2}(\lambda
y)-G_{2}(y)&=\beta\log(\lambda),\qquad \forall
\lambda>0.\end{split}
\end{equation}

\begin{equation}\label{eq22e:s}
\lim_{\theta\to 0}\widehat{G}_{2}(\lambda
\theta)-\widehat{G}_{2}(\theta)= -\beta\log(\lambda), \qquad
\forall \lambda>0.
\end{equation}
Both imply that
$$G_{2}(\theta)-\widehat{G}_{2}(1/\theta)\xrightarrow[\theta\to\infty]{}\beta \gamma.$$
\end{enumerate}
Where $\gamma$ is Euler's constant
$\gamma=\int^\infty_0e^{-v}\log(v)\mathrm{d}v.$
\end{lemma}

Tacking for granted Lemmas~\ref{lemmaML:s},\ref{lemmast:1} \&
\ref{lemmaAT:1} the proof of Theorem~\ref{teo1:s} is
straightforward.

A consequence of Lemma~\ref{lemmaAT:1} is that quantities related
to Sina\v\i 's condition can be used to determine whether or not
$\sigma$ has a finite expectation or a strictly positive drift.
That is the content of the following corollary.

\begin{cor}\label{cormeandrift}
\begin{enumerate}
\item[(i)]Assume that Sina\v\i 's condition holds at infinity with
index $\beta=1$ and that the lifetime of $\sigma$ is infinite.
Then $\sigma$ has a finite mean if and only if
$$R\equiv\lim_{\theta\to \infty}\log(\theta)-G_2(\theta)<\infty.$$ In
this case $\er(\sigma_1)=e^{\gamma+R}.$ \item[(ii)]Assume that
Sina\v\i 's condition holds at 0 with index $\beta=1.$ Then
$\sigma$ has a strictly positive drift $d$ if and only if
$$\widetilde{R}\equiv\lim_{\theta\to\infty}G_1(\theta)-\log(\theta)<\infty.$$ In
this case $d=e^{\gamma+\widetilde{R}}.$

\end{enumerate}
\end{cor}

\begin{remark}
For $\beta \in]0,1[,$ it is well known that $\phi\in
RV^\infty_\beta$ if and only if the sequence of subordinators
$\sigma^z$ defined by $(\sigma^z_t=z \sigma_{t/\phi(z)}, t\geq 0)$
converge as $z\to\infty$ in the sense of finite dimensional
distributions and in Skorohod's topology  to a stable subordinator
$\widetilde{\sigma}$ of parameter $\beta.$ This is equivalent to
say that for any $t>0$
\begin{equation}\label{remsinai1bis}
\pr(1<z\sigma_{t/\phi(z)}\leq
\lambda)\xrightarrow[z\to\infty]{}\pr(1<\widetilde{\sigma}_{t}\leq
\lambda),\qquad \lambda >1,
\end{equation} and
\begin{equation}\label{remsinai2bis}
\pr(\lambda<z\sigma_{t/\phi(z)}\leq
1)\xrightarrow[z\to\infty]{}\pr(\lambda<\widetilde{\sigma}_{t}\leq
 1),\qquad 0< \lambda<1.
\end{equation}
 On the other hand, Theorem~\ref{teo1:s} ensures that the
latter condition on $\phi$ holds if and only if Sina\v\i 's
condition holds at 0, which can be written as follows: for any
$\lambda>1$
\begin{equation}\label{remsinai1}
\begin{split}
\int^\infty_0\frac{\mathrm{d}s}{s}\pr(1<z^{-1}\sigma_{s/\phi(z^{-1})}\leq
\lambda )&=\int^\infty_0\frac{\mathrm{d}t}{t}\pr(z<\sigma_t\leq
\lambda z)
\\ &\xrightarrow[z\to 0]{} \beta \ln(\lambda)\\
&=\int^\infty_0\frac{\mathrm{d}s}{s}\pr(1<\widetilde{\sigma}_s\leq
\lambda),
\end{split}
\end{equation}
where the first equality is justified by a change of variables
$s=t\phi(z^{-1})$ and the last one follows from the scaling
property of the stable subordinator $\widetilde{\sigma};$ and for
any $0<\lambda<1$
\begin{equation}\label{remsinai2}
\begin{split}
\int^\infty_0\frac{\mathrm{d}s}{s}\pr(\lambda<z^{-1}\sigma_{s/\phi(z^{-1})}\leq
1)&\xrightarrow[z\to 0]{}
\int^\infty_0\frac{\mathrm{d}s}{s}\pr(\lambda<\widetilde{\sigma}_s\leq
1),
\end{split}
\end{equation}

Putting the pieces together we get that the result in
Theorem~\ref{teo1:s} can be viewed as an equivalence between the
convergence of the uni-dimensional laws of $\sigma^z$ in
(\ref{remsinai1bis}) \& (\ref{remsinai2bis}) and the convergence of
the integrated ones in (\ref{remsinai1}) \& (\ref{remsinai2}). An
analogous fact can be deduced for the convergence of $\sigma^z$ as
$z$ goes to 0 whenever Sina\v\i 's condition holds at infinity.
\end{remark}

We pass now to the proof of Lemmas~\ref{lemmaML:s},\ref{lemmast:1}
\& \ref{lemmaAT:1}.
\begin{subsection}{Proof of Lemmas~\ref{lemmaML:s},\ref{lemmast:1} \& \ref{lemmaAT:1}}
\begin{proof}[Proof of (i) in Lemma~\ref{lemmaML:s}]
We have by Frullani's formula that for every $\theta>0$
\begin{equation*}
\begin{split} 1+\phi(\theta)&=\exp\left\{\int^{\infty}_0\frac{\mathrm{d}t}{t}e^{- t}(1-e^{-t \phi(\theta)})\right\}\\
&= \exp\left\{\theta\int^{\infty}_0\mathrm{d}ye^{-\theta y}\int^{\infty}_0\frac{\mathrm{d}t}{t}e^{- t}\pr(\sigma_t>y)\right\}\\
&= \exp\{ \widehat{G}_{1}(\theta)\}.
\end{split}
\end{equation*}
\end{proof}

\begin{proof}[Proof of (ii) in Lemma~\ref{lemmaML:s}]
The equation relating $\phi$ and the measure $G_{2}$ can be
obtained using Frullani's formula but to prove moreover that this
measure is in fact is an harmonic renewal measure we proceed as
follows. Let $(\mathrm{e}_k, k\geq 1)$ be a sequence of
independent identically distributed random variables with
exponential law of parameter $1$  and independents of $\sigma.$
Put $\Theta_{l}=\sum^l_{k=1}\mathrm{e}_k, l\geq 1.$ It was proved
by Bertoin and Doney~\cite{MR96b:60190} that
$(\sigma_{\Theta_l},l\geq 1)$ forms a renewal process. The
harmonic renewal measure associated to $(\sigma_{\Theta_l},l\geq
1)$ is $G_{2}(\mathrm{d}x).$ Indeed,
\begin{equation*}
\begin{split}
\sum^{\infty}_{l=1}\frac{1}{l}\pr(\sigma_{\Theta_l}\in\mathrm{d}x)&=\sum^{\infty}_{l=1}\frac{1}{l}\int^{\infty}_0\mathrm{d}t \frac{t^{l-1}}{(l-1)!}e^{-t}\pr(\sigma_t \in\mathrm{d}x)\\
&=\int^{\infty}_0\frac{\mathrm{d}t}{t}e^{- t}(e^{
t}-1)\pr(\sigma_t\in\mathrm{d}x)=G_{2}(\mathrm{d}x).
\end{split}
\end{equation*}
Moreover, since the $l$-convolution of
$F(\mathrm{d}x)=\pr(\sigma_{\Theta_1}\in\mathrm{d}x)$ is such that
$$F^{*l}(\mathrm{d}x)=\pr(\sigma_{\Theta_l}\in\mathrm{d}x)$$ we
have that the Laplace transform $\widehat{F}(\theta)$ of $F$ is
related to that of $G_{2}$ by the formula
$$1-\widehat{F}(\theta)=\exp\{-\widehat{G}_{2}(\theta)\}\qquad
\theta>0.$$ Which finish the proof since $\widehat{F}(\theta)=
(1+\phi(\theta))^{-1},$ for $\theta>0.$ \end{proof}
\begin{proof}[Proof of (i) in Lemma \ref{lemmast:1}]
We can suppose without loss of generality that $\lambda>1.$ Given
that
\begin{equation*}
\begin{split}
\int^\infty_0\frac{\mathrm{d}t}{t}\pr(z < \sigma_t \leq \lambda z)= & \int^\infty_0\frac{\mathrm{d}t}{t}e^{- t}\pr(z < \sigma_t \leq \lambda z)+ \int^\infty_0\frac{\mathrm{d}t}{t}(1-e^{- t})\pr(z < \sigma_t \leq \lambda z),
\end{split}
\end{equation*}
for every $z>0,$ and
\begin{equation*}
\begin{split}
G_{1}(\lambda z)-G_{1}(z)&=\int^\infty_0\frac{\mathrm{d}t}{t}e^{- t}\pr(\frac{1}{\lambda z} \leq \sigma_t<\frac{1}{z}), \qquad z>0,
\end{split}
\end{equation*}
in order to prove (i) in Lemma~\ref{lemmast:1} we only need to check that
$$\lim_{z\to 0}\int^{\infty}_0\frac{\mathrm{d}t}{t}(1-e^{- t})\pr(z< \sigma_t\leq \lambda z)=0,\qquad \forall \lambda>1.$$
Indeed, given that for any $0<z<\infty$
$$\int^{1}_0\frac{\mathrm{d}t}{t}(1-e^{-t})\pr(\sigma_t\leq \lambda z)\leq \int^1_0\frac{\mathrm{d}t}{t}(1-e^{- t})<\infty,$$
we have by the monotone convergence theorem that
$$\int^{1}_0\frac{\mathrm{d}t}{t}(1-e^{- t})\pr(z < \sigma_t\leq \lambda z)\leq \int^{1}_0\frac{\mathrm{d}t}{t}(1-e^{- t})\pr(\sigma_t\leq\lambda z)\longrightarrow 0,\qquad\text{as} \ z \to 0, \quad \forall \lambda>1.$$
Furthermore, Theorem IV.7 in \cite{gihmanskorohod} enable us to
ensure that for any $0<z<\infty,$ and $1<\lambda,$
$$\int^{\infty}_1\frac{\mathrm{d}t}{t}\pr(\sigma_t \leq \lambda
z)<\infty.$$ Thus, proceeding as in the case $\int^1_0$ we obtain
that for any $\lambda >1,$
$$\lim_{z\to0}\int^{\infty}_1\frac{\mathrm{d}t}{t}(1-e^{- t})\pr(z<
\sigma_t\leq \lambda z)=0.$$
\end{proof}

\begin{proof}[Proof of (ii) in Lemma \ref{lemmast:1}] As in the proof of (i) it is enough to prove that
$$\lim_{z\to \infty}\int^{\infty}_0\frac{\mathrm{d}t}{t}e^{- t}\pr(z < \sigma_t\leq \lambda z)=0,\qquad \forall \lambda>1.$$
Indeed, it is straightforward that
$$\lim_{z\to\infty}\int^{\infty}_1\frac{\mathrm{d}t}{t}e^{- t}\pr(z
<\sigma_t\leq\lambda z)=0,\qquad \forall \lambda>1.$$ To prove
that
$$\lim_{z\to \infty}\int^{1}_0\frac{\mathrm{d}t}{t}e^{- t}\pr(z
<\sigma_t\leq\lambda z)=0,\qquad \forall \lambda>1.$$ we will use
the inequality~(6) in Lemma~1 of \cite{MR45:1250} which enable us
to ensure that for any $u>0$ and $z>0$
$$\pr(z <\sigma_t <\infty)\leq \frac{t \widetilde{\phi}(u) e^{-\kappa t}}{1-e^{-u
z}},\quad\text{with}\
\widetilde{\phi}(u)=du+\int_{]0,\infty[}(1-e^{-u x})
\pi(\mathrm{d}x).$$ Applying this inequality we get that, for any
$u,z>0,$ and $\lambda>1,$
$$\int^1_0\frac{\mathrm{d}t}{t}e^{-t}\pr(z <\sigma_t \leq \lambda
z) \leq \frac{\widetilde{\phi}(u)}{1-e^{-u z}}\int^1_0
e^{-(1+\kappa) t} \mathrm{d}t.$$ Making, first $z \to \infty$ and
then $u\to 0$  in the previous inequality, we obtain the estimate
$$\lim_{z\to\infty}\int^1_0\frac{\mathrm{d}t}{t}e^{-n t}\pr(z <
\sigma_t \leq \lambda z)\leq \widetilde{\phi}(0)\int^1_0
e^{-(1+\kappa) t} \mathrm{d}t,$$ valid for any $\lambda>1.$ Which
in fact ends the proof since
$\widetilde{\phi}(0)=\phi(0)-\kappa=0.$
\end{proof}

\begin{proof}[Proof of Lemma~\ref{lemmaAT:1}] The equivalence in (ii) of Lemma~\ref{lemmaAT:1} follows from Theorem 3.9.1 in~\cite{MR90i:26003}.
The equivalence in (i) of Lemma~\ref{lemmaAT:1} is obtained by applying Abelian-Tauberian theorems relying the behavior of the differences
of a non-decreasing function and those of their Mellin transform. Indeed, to prove that (\ref{eq2:s}) implies (\ref{eq2e:s}) we apply an Abelian
theorem that appears in \cite{MR90i:26003} Section 4.11.1. To that end we just need to verify that the Mellin transform $\check{k}$ of the kernel $k$ is
finite in a set $A=\{x\in\mathbb{C}: a\leq \Re(x)\leq b\}$ with $a<0<b.$ This is indeed the case since the Mellin transform of $k,$
$$\check k(x):=\int^\infty_0 t^{-x}k(t)\frac{\mathrm{d}t}{t}=\int^{\infty}_0t^{-x}e^{-t}\mathrm{d}t,\quad x\in\mathbb{C},$$ is finite in the strip $\Re(x)<1.$
That (\ref{eq2e:s}) implies (\ref{eq2:s}) is a direct consequence
of a Tauberian theorem for differences established
in~\cite{MR89a:26002}~Theorem~2.35.
\end{proof}

\begin{proof}[Proof of Corollary~\ref{cormeandrift}]
We will only prove the assertion in (i) of
Corollary~\ref{cormeandrift}. It is well known that any
subordinator has finite mean if and only if its Laplace exponent
is derivable at 0. Since $\sigma$ is assumed to have infinite
lifetime and the following relations, which are a consequence of
the Lemmas~\ref{lemmaML:s},\ref{lemmast:1} \& \ref{lemmaAT:1},
\begin{equation*}
\begin{split}
\frac{\phi(\theta)}{\theta} \sim
\frac{\displaystyle\frac{\phi(\theta)}{1+\phi(\theta)}}{\theta}
\sim \exp\left\{-\widehat{G}_2(\theta)-\log(\theta)\right\} \sim
\exp\left\{\gamma - G_2(1/\theta)+\log(1/\theta)\right\}, \quad
\text{as}\  \theta\to 0,
\end{split}
\end{equation*}
we have that $\phi$ is derivable at 0 if and only if the limit in
Corollary~\ref{cormeandrift}~(i) holds.

The proof of the assertion in Corollary~\ref{cormeandrift}~(ii)
uses the fact that $\sigma$ has a strictly positive drift if and
only if $\lim_{\theta\to\infty}\phi(\theta)/\theta>0$ and an
argument similar to that above. \end{proof}
\end{subsection}

\end{section}

\begin{section}{Proof of Theorem~\ref{teoch1:DA}}\label{proofthm2}
\begin{proof}[Proof of Theorem~\ref{teoch1:DA}]
This proof is a reworking of its analogous for random walks, which
was established by Rogozin~\cite{rogozin71} Theorem~9. We will
prove that under the assumptions of Theorem~\ref{teoch1:DA}, for
$t\to \infty$ in equation (\ref{eq3ch1:DA}), the assertion in
Corollary~\ref{corsec1:arcsin}~(iv) holds. (The proof of the case
$t\to 0$ in equation (\ref{eq3ch1:DA}) follows in a similar way
and so we omit the proof.) To that end, let $(\xi^r)_{r>0}$ be the
family of L\'evy processes defined by,
  $(\xi^r(t)=\xi_{rt}/b(r), t\geq 0)$ for $r>0.$ The hypothesis of
Theorem~\ref{teoch1:DA} is equivalent to the convergence, in the
sense of finite dimensional distributions, of the sequence of
L\'evy processes $\xi^r$ to a stable L\'evy process $X$ with
characteristic exponent given by the formula~(\ref{eq3ch1:DA}). By
Corollary~3.6 in Jacod--Shiryaev we have that this convergence
holds also in the Skorohod topology and Theorem~IV.2.3 in Gihman
\& Skorohod~\cite{gihmanskorohod} enable us to ensure that there
is also convergence of the first passage time above the level $x$
and the overshoot at first passage time above the level $x$ by
$\xi^r$ to the corresponding objects for $X.$ That is, for any
$x>0$
$$\tau^r_x=\inf\{t>0: \xi^r(t)>x\}, \qquad
\gamma^r_x=\xi^r(\tau^r_x)-x,$$ $$ \tau_x=\inf\{t>0:
X(t)>x\},\qquad  \gamma_x=X(\tau_x)-x,$$ we have that
$$(\tau^r_x,\gamma^r_x)\xrightarrow[r\to\infty]{\mathcal{D}}(\tau_x,\gamma_x).$$

In particular, for $x=1,$ we have that $\tau^r_1=r^{-1}T_{b(r)}$
and $\gamma^r_1=\left(\xi_{T_{b(r)}}-b(r)\right)/b(r),$ in the
notation of Corollary~\ref{corsec1:arcsin}, and thus that
$$\left(r^{-1}T_{b(r)},\left(\xi_{T_{b(r)}}-b(r)\right)/b(r)\right)\xrightarrow[r\to\infty]{\mathcal{D}}(\tau_1,\gamma_1).$$
We will next prove that
$$\left(\xi_{T_{r}}-r\right)/r\xrightarrow[r\to\infty]{\mathcal{D}}\gamma_1,$$
which implies that the assertion (iv) in
Corollary~\ref{corsec1:arcsin} holds. To that end, we introduce
the generalized inverse of $b,$ $b^{\leftarrow}(t)=\inf\{r>0:
b(r)>t\}$ for $t>0.$ Given that $b$ is regularly varying at
infinity it is known that $b(b^{\leftarrow}(t))\sim t$ as
$t\to\infty,$ see e.g.~\cite{MR90i:26003}~Theorem~1.5.12. Owing
the following relations valid for any $\epsilon>0$ fixed and small
enough,
$$b(b^{\leftarrow}(r)-\epsilon)\leq r \leq b(b^{\leftarrow}(r)),\qquad
r>0,$$we have that for any $x>0$
\begin{equation*}
\begin{split}
\pr\left(\frac{\xi_{T_{b(b^{\leftarrow}(r))}}}{b(b^{\leftarrow}(r))}\frac{b(b^{\leftarrow}(r))}{r}\leq
x+1\right)& \leq \pr\left(\frac{\xi_{T_{r}}}{r}\leq
x+1\right)\\
&\leq
\pr\left(\frac{\xi_{T_{b(b^{\leftarrow}(r)-\epsilon)}}}{b(b^{\leftarrow}(r)-\epsilon)}\frac{b(b^{\leftarrow}(r)-\epsilon)}{b(b^{\leftarrow}(r))}\frac{b(b^{\leftarrow}(r))}{r}\leq
x+1\right).
\end{split}
\end{equation*}
Making $r$ tend to infinity and using that
$b(b^{\leftarrow}(r)-\epsilon)/b(b^{\leftarrow}(r))\to 1$ as $r\to
\infty,$ we get that the left and right hand sides of the previous
inequality tend to $\pr(\gamma_1+1\leq x+1)$ and so that for any
$x>0$
$$\pr\left(\frac{\xi_{T_{r}}-r}{r}\leq
x\right)\xrightarrow[r\to\infty]{}\pr(\gamma_1\leq x).$$

Furthermore, it is well known that in the case
$\alpha\rho\in]0,1[$ the law of $\gamma_1$ is the generalized
arc-sine law with parameter $\alpha\pr(X_1>0)=\alpha\rho,$ that is
$$\pr(\gamma_1\in \mathrm{d}x)=\frac{\sin(\alpha\rho
\pi)}{\pi}x^{-\alpha\rho}(1+x)^{\alpha\rho-1}\mathrm{d}x,\quad
x>0.$$ In the case $\alpha\rho=0$ the random variable $\gamma_1$
is degenerate at infinity and  in the case $\alpha\rho=1$ it is
degenerate at $0$. Thus, in any case the Sinai index of $\xi$ is
$\alpha\rho.$
\end{proof}

For shake of completeness in the following Lemma we provide
necessary and sufficient conditions on the tail behavior of the
L\'evy measure of $\xi$ in order that the hypotheses of
Theorem~\ref{teoch1:DA} be satisfied. This result concerns only
the case $t\to\infty$ in (i) of Theorem~\ref{teoch1:DA} and
$\alpha\in]0,1[.$ The triple $(a,q^2,\Pi)$ denotes the
characteristics of the L\'evy process $\xi,$ that is, its linear
and Gaussian term, $a,q$ and L\'evy measure $\Pi$ and are such
that $$\Psi(\lambda)=ia\lambda-\frac{\lambda^2
q^2}{2}+\int_{\re\setminus\{0\}}(e^{i\lambda x}-1-i\lambda
x1_{\{\vert x\vert<1\}})\Pi(\mathrm{d}x),\qquad \lambda\in\re.$$
By $\overline{\Pi}^+$ and $\overline{\Pi}^-$ we denote the right
and left hand tails of the L\'evy measure $\Pi$ respectively, i.e.
$\overline{\Pi}^+(x)=\Pi]x,\infty[$ and
$\overline{\Pi}^-(x)=\Pi]-\infty,-x[,$ for $x>0.$

\begin{lemma}\label{lemDA} Let $\alpha\in]0,1[$ and $\delta\in[-1,1].$ The following are equivalent
\begin{description}
\item[DA] There exits a function $b:]0,\infty[\to]0,\infty[$ which
is regularly varying at infinity with index $\beta=1/\alpha$ and
such that the limit in equation~(\ref{eq3ch1:DA}) holds as
$t\to\infty$. \item[TB] The function $\overline{\Pi}^+(\cdot) +
\overline{\Pi}^-(\cdot)$ is regularly varying at infinity with
index $-\alpha$ and
\begin{equation*}
\frac{\overline{\Pi}^+(x)}{\overline{\Pi}^+(x) + \overline{\Pi}^-(x)}\longrightarrow p, \quad \frac{\overline{\Pi}^-(x)}{\overline{\Pi}^+(x) + \overline{\Pi}^-(x)}\longrightarrow q, \quad \text{as}\quad x \to \infty; \quad p+q=1, p-q=\delta.
\end{equation*}
\end{description}
\end{lemma}

\begin{remark}The same result holds true for $\alpha\in]0,2[$ if the L\'evy process is assumed to be symmetric (the proof of Lemma~\ref{lemDA} can be easily extended to this case). Furthermore, there is also an analogue of this result when $t\to 0$ in (i) of Theorem~\ref{teoch1:DA} in the cases $1<\alpha<2$ or $0<\alpha<2$ and $\xi$ is assumed to be symmetric. Its  proof is quite similar to that of Lemma~\ref{lemDA}, see e.g. the recent work of De Weert~\cite{Deweert}.
\end{remark}

\begin{proof}[Proof of Lemma~\ref{lemDA}]
It is plain, that for any $t>0$ the function
$\Psi^{(t)}(\lambda):=
t\Psi\left(\displaystyle{\frac{\lambda}{b(t)}}\right)$ is the
characteristic exponent of the infinitely divisible random
variable $X^{(t)}:=\xi_t/b(t),$ which by the hypothesis
DA($\alpha$) converges to a stable law $X(1)$ whose characteristic
exponent is given by equation~(\ref{eq3ch1:DA}). The
characteristic exponent $\Psi^{(t)}$ can be written as

\begin{equation*}
\Psi^{(t)}(\lambda)=i\lambda
a^{(t)}-\lambda^2(q^{(t)})^2/2+\int_{\re\setminus\{0\}}\left(e^{i\lambda
z}-1-ih(z)\right)\Pi^{(t)}(\mathrm{d}z),
\end{equation*} where $h(z)=z1_{\{|z|\leq 1\}}+z^{-1}1_{\{|z|> 1\}}$ and $(a^{(t)}, q^{(t)}, \Pi^{(t)})$ are given by

\begin{gather*} a^{(t)}=\frac{t a}{b(t)}+\frac{t}{b(t)}\int x1_{\{1< |x|\leq b(t)\}} \Pi(\mathrm{d}x)+ tb(t)\int x^{-1}1_{\{|x|>b(t)\}}\Pi(\mathrm{d}x),\\
q^{(t)}=q \left(\frac{t}{b(t)^2}\right)^{1/2},\quad
\Pi^{(t)}(\mathrm{d}z)=t\Pi(b(t)\mathrm{d}z).
\end{gather*}

According to a well known result on the convergence of infinite
divisible laws, see e.g. Sato~\cite{sato99}, the convergence in
law of $X^{(t)}$ to $X$ as $t\to\infty,$ is equivalent to the
convergence of the triplet $(a^{(t)}, q^{(t)}, \Pi^{(t)})$ to
$(l,0,\Pi_S)$ as $t\to \infty,$ with
$$\Pi_S(\mathrm{d}x)=\left(c_+x^{-1-\alpha}1_{\{x>0\}}+c_-|x|^{-1-\alpha}1_{\{x<0\}}\right)\mathrm{d}x,\qquad
c_+, c_-\in\re^+,$$ and $$l=\frac{2(c_+ - c_-)}{1-\alpha^2},$$ and
they are such that for every $\lambda\in\re$
\begin{equation*}\begin{split}-c|\lambda|^\alpha\left(1-i\delta\textup{sgn}(\lambda)\tan(\pi\alpha/2)\right)&=\int_{\re\setminus\{0\}}(e^{i\lambda x}-1)\Pi_S(\mathrm{d}x)\\
&=il+\int_{\re\setminus\{0\}}\left(e^{i\lambda
x}-1-ih(x)\right)\Pi_S(\mathrm{d}x),\end{split} \end{equation*}
with $c>0$ convenably chosen. The term $\delta$ in the previous
equation is determined by $\delta=p-q=\displaystyle{\frac{c_+ -
c_-}{c_+ + c_-}}.$

That the hypotheses on the tail behavior of $\Pi$ are equivalent
to the convergence of $\Pi^t$ to $\Pi_S$ is a quite standard fact
in the theory of domains of attraction and so we refer
to~\cite{MR90i:26003} section 8.3.2, for a proof. This implies in
particular that for any $x>0,$
$$t\overline{\Pi}^{+}(b(t)x)\rightarrow
c_+x^{-\alpha},\quad\text{and}\quad
t\overline{\Pi}^{-}(b(t)x)\rightarrow
c_-x^{-\alpha},\quad\text{as}\  t\to\infty.$$ The only technical
detail that requires a proof is that $a^{(t)}\to l$ as
$t\to\infty.$ Indeed, under the conditions (ii) of
Theorem~\ref{teoch1:DA} and $0<p<1$ the functions
$\overline{\Pi}^+(\cdot)$ and $\overline{\Pi}^-(\cdot)$ are
regularly varying at infinity with index $0<\alpha<1,$ this
implies that $a t/b(t)\to 0$ as $t\to\infty.$ Moreover, it is
justified by making an integration by parts that
\begin{equation*}
\begin{split}
\int x1_{\{1< x\leq b(t)\}}\Pi(\mathrm{d}x)&\sim \overline{\Pi}^+(1)-b(t)\overline{\Pi}^+(b(t))+\int^{b(t)}_1\overline{\Pi}^+(z)dz\\
&\sim
\overline{\Pi}^+(1)-b(t)\overline{\Pi}^+(b(t))+\frac{b(t)\overline{\Pi}^+(b(t))}{1-\alpha},
\end{split}
\end{equation*} as $t\to\infty.$ Multiplying by $t/b(t)$ we get that
\begin{equation*}
\begin{split}
\frac{t}{b(t)}\int x1_{\{1< x\leq b(t)\}}\Pi(\mathrm{d}x)&\sim
-t\overline{\Pi}^+(b(t))+\frac{t\overline{\Pi}^+(b(t))}{1-\alpha}\longrightarrow
-c_+ + \frac{c_+}{1-\alpha}
\end{split}
\end{equation*} as $t\to\infty.$ Similarly, it is proved that $$\frac{t}{b(t)}\int x1_{\{1< -x\leq b(t)\}}\Pi(\mathrm{d}x)\vat c_- - \frac{c_-}{1-\alpha}$$
Concerning the term $\int x^{-1}1_{\{|x|>b(t)\}}\Pi(\mathrm{d}x),$
an integration by parts and Karamata's theorem yield
\begin{equation*}
\begin{split}
\int x^{-1}1_{\{x> b(t)\}}\Pi(\mathrm{d}x)&\sim (b(t))^{-1}\overline{\Pi}^+(b(t))+\int^{\infty}_{b(t)}z^{-2}\overline{\Pi}^+(z)dz\\
&\sim
(b(t))^{-1}\overline{\Pi}^+(b(t))+\frac{\overline{\Pi}^+(b(t))}{b(t)(1+\alpha)},
\end{split}
\end{equation*}
as $t\to\infty$ and therefore
\begin{equation*}
\begin{split}
t b(t)\int x^{-1}1_{\{x> b(t)\}}\Pi(\mathrm{d}x) &\sim
t\overline{\Pi}^+(b(t))+\frac{t\overline{\Pi}^+(b(t))}{1+\alpha}\longrightarrow
c_+ + \frac{c_+}{1+\alpha}.
\end{split}
\end{equation*} Analogously, we prove
\begin{equation*}
\begin{split}
t b(t)\int x^{-1}1_{\{x< -b(t)\}}\Pi(\mathrm{d}x)\van -c_- -
\frac{c_-}{1+\alpha}.
\end{split}
\end{equation*}Finally, adding up these four terms it follows that
$$\lim_{t\to\infty}a^{(t)}=\frac{c_+-c_-}{1-\alpha}+ \frac{c_+-c_-}{1+\alpha}=l.$$ The proof that $a^{(t)}\to l$ in the case $p=1,$ respectively $p=0,$ is quite similar but uses that $\overline{\Pi}^-=o(\overline{\Pi}^+),$ respectively $\overline{\Pi}^+=o(\overline{\Pi}^-).$
\end{proof}
\end{section}

\begin{section}{Some criteria based on the characteristics of
$\xi$}\label{somecriterions}

The result obtained in Theorem~\ref{mainthm:1} gives a complete
solution to the question posed at the introduction of this note.
Nevertheless, this answer is not completely satisfactory since the
criterion provided by this theorem is not given in terms of the
characteristics of $\xi.$  However, to get a result based on the
characteristics of $\xi$ one should be able to control the
behavior of the dual ladder height subordinator $\widehat{H}$,
that is, the ladder height subordinator of the dual L\'evy process
$\widehat{\xi}=-\xi.$ This is due to the fact, showed by
Vigon~\cite{MR2002i:60101}, that the L\'evy measure of the ladder
height subordinator $H$ is determined by the L\'evy measure of
$\xi$ and the potential measure of $\widehat{H}.$ (See the
Lemma~\ref{vigon} below for a precise statement.) With this in
mind we will make some assumptions on the dual ladder height
process, that can be verified directly from the characteristics of
$\xi,$ and provide some nasc for the regular variation, at
infinity or 0, of the Laplace exponent of the ladder height
subordinator $H$. To that end we will start by obtaining some tail
estimates of the right tail of the L\'evy measure of $H$ which
will provide us the necessary tools to determine whether or not
the underlying L\'evy process $\xi$ satisfies Sina\v\i 's
condition. To tackle this task we will first introduce some
notation and recall some known facts.

We will denote by $(k_0,d,po)$ (resp.
$(\widehat{k}_0,\widehat{d},ne)$) the characteristics of the
subordinator $H$ (resp. $\widehat{H}$) that is, its killing term,
drift and L\'evy measure, respectively. Let $\widehat{V}$ be the
potential measure of $\widehat{H},$ that is
$\widehat{V}(\mathrm{d}x)=\er(\int^{\infty}_0
1_{\{\widehat{H}_t\in \mathrm{d}x\}}\mathrm{d}t).$ Furthermore,
throughout this section we will denote by $(a,q,\Pi)$ the
characteristics of the L\'evy process $\xi.$ Finally, by the
symbols $\overline{po},\overline{ne},\overline{\Pi}^+,$ we denote
the right tail of the L\'evy measures of $H,\widehat{H}$ and $\xi$
respectively, that is $$\overline{po}(x)=po]x,\infty[,\quad
\overline{ne}(x)=ne]x,\infty[,\qquad
\overline{\Pi}^+(x)=\Pi]x,\infty[,\qquad x>0,$$ and by $\Pi^+$ the
restriction of $\Pi$ to $]0,\infty[,$
$\Pi^+=\Pi\vert_{]0,\infty[}$.

We recall that the Laplace exponent, $\varphi,$ of $H$ varies
regularly at 0 with index $\alpha\in]0,1[$ if and only if the
function $\overline{po}$ is regularly varying at infinity with
index $-\alpha.$ Analogously, $\varphi$ varies regularly at
infinity with index $\alpha\in]0,1[$ if and only if the drift
coefficient $d$ is zero and  the function $\overline{po}$ is
regularly varying at 0 with index $-\alpha.$ Owing this relations
we will restrict ourselves to study the behavior of the function
$\overline{po}.$

As we mentioned before, Vigon~\cite{MR2002i:60101} established
some identities ``equations amicales'' relying the L\'evy measures
$po, ne$ and $\Pi;$ these are quoted below for ease of reference.

\begin{lemma}[Vigon~\cite{MR2002i:60101}, Equations amicales]\label{vigon}
We have the following relations

\begin{enumerate}
\item[(EAI)] $\displaystyle
\overline{po}(x)=\int^{\infty}_0\widehat{V}(\mathrm{d}y)\overline{\Pi}(x+y),\qquad
x>0.$

\item[(EA)] $\displaystyle
\overline{\Pi}^+(x)=\int_{]x,\infty[}po(\mathrm{d}y)\overline{ne}(y-x)+\widehat{d}\overline{p}(x)+\widehat{k}_0\overline{po}(x)$,
for any $x>0;$ where $\overline{p}(x)$ is the density of the
measure $po,$ which exists if $\widehat{d}>0.$
\end{enumerate}
\end{lemma}

We will assume throughout that the underlying L\'evy process is
not spectrally negative, that is $\Pi]0,\infty[>0,$ since in that
case the ladder height process $H$ is simply a drift, $H_t=d t,
t\geq 0.$

We will say that a measure $\mathit{M}$ on $[0,\infty[$ belongs to
the class $\mathcal{L}^0$ of \textit{long tailed} measures if its
tail $\overline{\mathit{M}}(x)=\mathit{M}]x,\infty[,$ is such that
$0<\overline{\mathit{M}}(x)<\infty$ for each $x>0$ and
$$\lim_{x\to\infty}\frac{\overline{\mathit{M}}(x+t)}{\overline{\mathit{M}}(x)}=1,
\qquad \textit{for each}\ t\in\re.$$ It is well known that this
family includes the subexponential measures and the cases at which $\overline{M}$ is regularly varying.

In the first results of this section we relate the behavior of
$\overline{\Pi}^+$ with that of $\overline{po}$ at infinity. These
results are close in spirit to those obtained by Kl\"uppelberg,
Kyprianou and Maller~\cite{KKM}. In that work the authors assume
that the L\'evy process $\xi$ drifts to $-\infty,$ i.e.
$\lim_{t\to\infty}\xi_t=-\infty,$ $\pr$--a.s. and obtain several
asymptotic estimates of the function $\overline{po}$ in terms of
the L\'evy measure $\Pi.$ In our setting we permit any  behavior
of $\xi$ at the price of making some assumptions on the dual
ladder height subordinator. The following results are the
continuous time analogue of the result of
Veraverbeke~\cite{MR0423543} and Gr\"ubel~\cite{Grubel85} for
random walks.

\begin{teo}\label{chnasc:1}
\begin{enumerate}
\item[(a)]Assume that the dual ladder height subordinator has a finite mean $\mu=\er(\widehat{H}_1)<\infty.$ The following are equivalent
\begin{itemize}
\item[(a-1)]The measure $\Pi^+_{I}$ on $[0,\infty[$ with tail $\overline{\Pi}^+_{I}(x)=\int^{\infty}_{x}\overline{\Pi}^+(z)\mathrm{d}z, x\geq 0,$ belongs to $\mathcal{L}^0.$
\item[(a-2)] $po\in \mathcal{L}^0.$
\item[(a-3)] ${\displaystyle \overline{po}(x)\sim \frac{1}{\mu} \overline{\Pi}^+_{I}(x),}$ as $x\to \infty.$
\end{itemize}
\item[(b)] Assume that the dual ladder height subordinator $\widehat{H}$ has killing term $\widehat{k}_0 > 0.$ The following are equivalent
\begin{itemize}
\item[(b-1)] $\Pi^+=\Pi |_{]0,\infty[}\in \mathcal{L}^0.$
\item[(b-2)] $po\in\mathcal{L}^0,$ \begin{equation}\label{taub:s3}
\frac{\widehat{d}\overline{p}(x)}{\overline{po}(x)}\longrightarrow
0 \quad\text{and}\quad
\int^{1}_{0}\left(\frac{\overline{po}(x)-\overline{po}(x+y)}{\overline{po}(x)}\right)ne(\mathrm{d}y)\longrightarrow
0,\quad \text{as}\ x\to \infty.
\end{equation}
\item[(b-3)] ${\displaystyle \overline{po}(x)\sim \frac{1}{\widehat{k}_0} \overline{\Pi}^+(x)},$ as $x\to\infty.$
\end{itemize}
\end{enumerate}
\end{teo}
As a corollary of the previous Theorem we have the following
criterions to determine whether or not the tail of the L\'evy
measure of $H$ is regularly varying.
\begin{cor}\label{regvar:cr1}
\begin{enumerate}
\item[(a)]Under the assumptions of (a) in Theorem~\ref{chnasc:1}
and for any $\alpha\in]0,1]$ we have that $\overline{\Pi}^+\in
RV^{\infty}_{-1-\alpha}$ if and only if $\overline{po}\in
RV^{\infty}_{-\alpha}.$ Both imply that $$\overline{po}(x)\sim
\frac{1}{\alpha\mu}x\overline{\Pi}^+(x),\qquad x\to \infty.$$
\item[(b)]Under the assumptions of (b) in Theorem~\ref{chnasc:1}
and for any $\alpha\in]0,1[$ we have that $\overline{\Pi}^+\in
RV^{\infty}_{-\alpha}$ if and only if $\overline{po}\in
RV^{\infty}_{-\alpha}$ and \begin{equation*}
\int^{1}_{0}\left(\frac{\overline{po}(x)-\overline{po}(x+y)}{\overline{po}(x)}\right)ne(\mathrm{d}y)\longrightarrow
0,\quad \text{as}\ x\to \infty.
\end{equation*}
\end{enumerate}
\end{cor}
\begin{proof}
The proof of (a) in Corollary~\ref{regvar:cr1} follows from the
fact that under these hypotheses $$\frac{x
\overline{\Pi}^+(x)}{\int^\infty_x\overline{\Pi}^+(z)\mathrm{d}z}\longrightarrow
\alpha, \qquad \text{as}\ x\to \infty,$$ which is a consequence of
Theorem 1.5.11 in~\cite{MR90i:26003}. The proof of (b) is
straightforward.
\end{proof}

The behavior at 0 of $\overline{po}$ was studied by Vigon
in~\cite{thesevigon} Theorems 6.3.1 and 6.3.2. He obtained several
estimations that we will use here to provide an analogue of
Theorem~\ref{chnasc:1} for the behavior at 0 of $\overline{po}.$
(Those estimates are more general than needed, see the proof of
Proposition~\ref{bhvat0}.)

\begin{prop}\label{bhvat0}
\begin{enumerate}
\item[(a)] Assume that $\widehat{H}$ has a drift $\widehat{d}>0.$
For any $\alpha\in]0,1]$ we have that  $\overline{po}\in
RV^{0}_{-\alpha}$ if and only if $\overline{\Pi}^+\in
RV^{0}_{-\alpha-1}.$ \item[(b)] Assume that $\widehat{H}$ has a
drift $\widehat{d}=0$ and that the total mass of the measure $ne$
is finite, equivalently, $\lim_{x\to 0+}\overline{ne}(x)<\infty.$
Then for any $\alpha\in]0,1]$ we have that $\overline{po}\in
RV^{0}_{-\alpha}$ if and only if $\overline{\Pi}^+\in
RV^{0}_{-\alpha}.$ Moreover, the same assertion holds if
furthermore $\alpha=0$ and $\lim_{x\to0+}\overline{po}(x)=\infty.$
\end{enumerate}
\end{prop}

Before passing to the proof of Theorem~\ref{chnasc:1} and
Proposition~\ref{bhvat0} we will make some remarks.

\begin{remark} The assumptions in Theorem~\ref{chnasc:1} can be verified using only the characteristics of the underlying L\'evy process $\xi.$ According to Doney and Maller~\cite{MR2003d:60090}, necessary and sufficient conditions on $\xi$ to be such that $\er(\widehat{H}_1)<\infty,$ are either $0<\er(-\xi_1)\leq \er|\xi_1|<\infty$ or $0=\er(-\xi_1)<\er|\xi_1|<\infty$ and
$$\int_{[1,\infty[}\left(\frac{x\overline{\Pi}^-(x)}{1+\int^x_0\mathrm{d}y\int^{\infty}_y\overline{\Pi}^+(z)\mathrm{d}z}\right)\mathrm{d}x<\infty \quad \text{with}\ \overline{\Pi}^-(x)=\Pi]-\infty,-x[,\ x>0.$$
Observe that under such assumptions the L\'evy process $\xi$ does not drift to $\infty,$ i.e. $\liminf_{t\to\infty}\xi_t=-\infty,$ $\pr$--a.s. The case where the L\'evy process $\xi$ drift to $\infty,$ $\lim_{t\to\infty}\xi_t=\infty,$ $\pr$--a.s. or equivalently $\widehat{k}_0>0$ is considered in (b). It was proved by Doney and Maller that $\xi$ drift to $\infty$ if and only if
$$\int_{]-\infty,-1[}\left(\frac{|y|}{\overline{\Pi}^+(1)+\int^{|y|}_1\overline{\Pi}^+(z)\mathrm{d}z}\right)\Pi(\mathrm{d}y)<\infty=\int^{\infty}_1\overline{\Pi}^+(x)\mathrm{d}x \quad\text{or}\quad 0<\er(\xi_1)\leq \er|\xi_1|<\infty.$$

The assumptions in Proposition~\ref{bhvat0} can be verified using the recent results of Vigon~\cite{MR2002i:60101} and \cite{thesevigon} Chapter~10.
\end{remark}

\begin{remark}
The results in Theorem~\ref{chnasc:1} concern only the case at
which the underlying L\'evy process does not has exponential
moments and so it extends to L\'evy processes the Theorem~1-(B,C)
of Veraverbeke~\cite{MR0423543}. The case at which the L\'evy
process has exponential moments has been considered by
Kl\"uppelberg et al.~\cite{KKM} Proposition~5.3 under the
assumption that the underlying L\'evy process has positive jumps
and drifts to $-\infty,$ but actually the latter hypothesis is not
used in their proof, and so their result is still true in this
more general setting, which extends Theorem~1-A
in~\cite{MR0423543}.
\end{remark}

\begin{remark} The estimate of $\overline{po}$ obtained in Theorem~\ref{chnasc:1}-a holds whenever the function $\overline{\Pi}^+_I$ belongs to the class $\mathcal{L}^0,$ but it is known that this can occurs even if $\Pi^+\notin\mathcal{L}^0$, see e.g. Kl\"uppelberg~\cite{kluppelberg88}. A question arises: Is it possible to sharpen the estimate of $\overline{po}$ provided in Theorem~\ref{chnasc:1}-(a) when moreover $\Pi^+\in\mathcal{L}^0$? The following result answers this question in affirmatively.

\begin{prop}\label{keythm}
Assume that $\mu=\er(\widehat{H}_1)<\infty.$ The following are equivalent
\begin{enumerate}
\item[(i)] $\displaystyle \Pi^+\in \mathcal{L}^0$.
\item[(ii)] For any $g:\re^+\to\re^+$ directly Riemman integrable,
$$\displaystyle \lim_{x\to\infty}\frac{1}{\overline{\Pi}^+(x)}\int^{\infty}_x po(\mathrm{d}y)g(y-x) = \frac{1}{\mu}\int^{\infty}_{0}g(z)\mathrm{d}z.$$
\end{enumerate}
\end{prop}

To our knowledge the discrete time analogue of this result, that
we state below, is unknown, although it can be easily deduced from
the arguments in Asmussen et al.~\cite{asmussenetal2002} Lemma~3.
We use the notation in the introduction of this work.

\begin{prop}\label{keythmdisc}
Assume that $m=\er(Z_{\widehat{N}})<\infty,$ where
$\widehat{N}=\inf\{n>0: Z_n\leq 0\}$ and that the law of $X_1$ is
non-lattice.  The following are equivalent
\begin{enumerate}
\item[(i)] The law of $X_1$ belongs to the class $\mathcal{L}^0$.
\item[(ii)] For any $g:\re^+\to\re^+$ directly Riemman integrable,
$$\displaystyle \lim_{x\to\infty}\frac{1}{\overline{F}(x)}\int^{\infty}_x g(y-x) \pr(Z_{N}\in\mathrm{d}y) = \frac{1}{m}\int^{\infty}_{0}g(z)\mathrm{d}z,$$
where $\overline{F}(x)=\pr(X_1>x),\ x>0.$
\end{enumerate}
\end{prop}
The proof of this result is quite similar to that of
Proposition~\ref{keythm} and so we will omit it.
\end{remark}

\begin{subsection}{Proofs of Theorem~\ref{chnasc:1} and Propositions~\ref{bhvat0} and~\ref{keythm}}

\begin{proof}[Proof of (a) in Theorem~\ref{chnasc:1}] To prove that (a-1) is equivalent to (a-2) we will prove that either of this conditions implies that
\begin{equation}\label{eqpo1:1}
\overline{po}(x)=\int^\infty_0\overline{\Pi}^+(x+y)\widehat{V}(\mathrm{d}y)\sim \frac{1}{\mu}\int^\infty_x\overline{\Pi}^+(z)\mathrm{d}z, \qquad \text{as} \ x\to\infty,
\end{equation} with $\mu:=\er(\widehat{H}_1),$ from where the result follows. (Observe that the assumption that $\widehat{H}$ has a finite mean implies that $\int^\infty\overline{\Pi}^+(z)\mathrm{d}z<\infty.$)

Assume that (a-1) holds. Indeed,  by the renewal theorem for
subordinators we have that for any $h>0,$
$$\lim_{t\to\infty}\widehat{V}]t,t+h]=\frac{h}{\mu}.$$ Thus, for
any $h>0$ given and any $\epsilon>0$ there exists a
$t_0(h,\epsilon)>0$ such that
$$(1-\epsilon)\frac{h}{\mu}<\widehat{V}]t,t+h]<(1+\epsilon)\frac{h}{\mu},\qquad
\forall t>t_0,$$ and as a consequence, if $N_0$ is an integer such
that $N_0h>t_0,$ we have the following inequalities
\begin{equation*}
\begin{split}
\int^\infty_0\overline{\Pi}^+(x+y)\widehat{V}(\mathrm{d}y)&\leq \overline{\Pi}^+(x)\widehat{V}[0,N_0h]+\sum^{\infty}_{k=N_0} \overline{\Pi}^+(kh + x) \widehat{V}]kh,kh+h]\\
&\leq \overline{\Pi}^+(x)\widehat{V}[0,N_0h]+(1+\epsilon) \sum^{\infty}_{k=N_0} \overline{\Pi}^+(kh + x)\frac{h}{\mu}\\
&\leq \overline{\Pi}^+(x)\widehat{V}[0,N_0h]+\frac{(1+\epsilon)}{\mu}\int^{\infty}_{(N_0-1)h}\overline{\Pi}^+(x+z)\mathrm{d}z\\
&\leq \overline{\Pi}^+(x)\widehat{V}[0,N_0h]+\frac{(1+\epsilon)}{\mu}\int^{\infty}_{x}\overline{\Pi}^+(z)\mathrm{d}z.
\end{split}
\end{equation*}

It follows from the previous inequalities and the fact that
$$\overline{\Pi}^+(x)/\int^{\infty}_{x}\overline{\Pi}^+(z)\mathrm{d}z\longrightarrow
0,\qquad \text{as} \ x \to \infty,$$ since $\overline\Pi^+_I\in
\mathcal{L}^0$ and $\overline{\Pi}^+$ is decreasing, that
$$\limsup_{x\to\infty}
\frac{\int^\infty_0\overline{\Pi}^+(x+y)\widehat{V}(\mathrm{d}y)}{\frac{1}{\mu}\int^{\infty}_{x}\overline{\Pi}^+(z)\mathrm{d}z}\leq
1.$$ Analogously, we prove that
$$\int^\infty_0\overline{\Pi}^+(x+y)\widehat{V}(\mathrm{d}y)\geq
\frac{(1-\epsilon)}{\mu}\int^\infty_x\overline{\Pi}^+(z)\mathrm{d}z
- \frac{(1-\epsilon)}{\mu}\overline{\Pi}^+(x) (N_0+1)h,  \qquad
x>0.$$ Therefore,
$$\liminf_{x\to\infty} \frac{\int^\infty_0\overline{\Pi}^+(x+y)\widehat{V}(\mathrm{d}y)}{\frac{1}{\mu}\int^{\infty}_{x}\overline{\Pi}^+(z)\mathrm{d}z}\geq 1.$$ Which ends the proof of the claim~(\ref{eqpo1:1}).

We assume now that (a-2) holds and we will prove that the estimate
in~(\ref{eqpo1:1}) holds. On the one hand, we know that for every
$z>0$ $$\overline{\Pi}^+(z)=\int^\infty_z
po(\mathrm{d}y)\overline{ne}(y-z)+\widehat{d}\overline{p}(z),$$
since $\widehat{k}_0=0,$ because under our assumptions the L\'evy
process does not drift to $\infty. $ Integrating this relation
between $x$ and $\infty$ and using Fubini's theorem we obtain that
for any $x>0$
\begin{equation*}
\begin{split}
\int^\infty_x \mathrm{d}z \overline{\Pi}^+(z)&=\int^\infty_x po(\mathrm{d}y)\int^{y-x}_0 \mathrm{d}z \overline{ne}(z) + \widehat{d}\overline{po}(x)\\
&\leq \overline{po}(x)\left(\int^{\infty}_0\mathrm{d}z\overline{ne}(z)+\widehat{d}\right)\\
&=\overline{po}(x)\mu<\infty.
\end{split}
\end{equation*}
Thus, $$\limsup_{x\to\infty}\frac{\frac{1}{\mu}\int^\infty_x
\mathrm{d}z \overline{\Pi}^+(z)}{\overline{po}(x)}\leq 1.$$ On the
other hand, to prove that
$$\liminf_{x\to\infty}\frac{\frac{1}{\mu}\int^\infty_x \mathrm{d}z \overline{\Pi}^+(z)}{\overline{po}(x)}\geq 1,$$
we will use an argument based on some facts of renewal theory.
To that end we recall that it was proved by Bertoin and Doney~\cite{MR96b:60190} that the potential
measure of a subordinator $\widehat{H}$ is the delayed renewal measure associated to the law
$F(x)=\pr(\widehat{H}_\vartheta\leq x)$ with $\vartheta$ an exponential random variable independent of $\widehat{H},$
that is
$$\widehat{V}(\mathrm{d}y)=\sum^\infty_{n=1}F^{*n}(\mathrm{d}y).$$
We have by hypothesis that $\int^\infty_0(1-F(x))\mathrm{d}x=\er(\widehat{H}_1)=\mu<\infty$
and thus the measure $\widetilde{G}_F(\mathrm{d}y)$ on $]0,\infty[,$ with density
$G_F(z):= (1-F(z))/\mu, z>0$ is a probability measure. By standard facts of renewal theory we know that the following equality between measures holds
$$\frac{\mathrm{d}y}{\mu}=\widetilde{G}_F(\mathrm{d}y)+\widetilde{G}_F*\widehat{V}(\mathrm{d}y),\qquad y>0,$$
where $*$ denotes the standard convolution between measures. Using this identity and the equation (EAI) we have that for any $x>0,$
\begin{equation*}
\begin{split}
\frac{1}{\mu}\int^\infty_0\mathrm{d}y \overline{\Pi}^+(x+y)=&\int^\infty_0\mathrm{d}y G_F(y)\overline{\Pi}^+(x+y)+ \int^\infty_0\mathrm{d}z G_F(z)\int^\infty_0\widehat{V}(dr)\overline{\Pi}^+(x+z+r)\\
=&\int^\infty_0\mathrm{d}y G_F(y)\overline{\Pi}^+(x+y)+ \int^\infty_0\mathrm{d}z G_F(z)\overline{po}(x+z),
\end{split}
\end{equation*} and by Fatou's lemma we get that $$\liminf_{x\to\infty}\frac{\frac{1}{\mu}\int^\infty_0\mathrm{d}y \overline{\Pi}^+(x+y)}{\overline{po}(x)}\geq \int^\infty_0\mathrm{d}z G_F(z)\liminf_{x\to\infty}\frac{\overline{po}(x+z)}{\overline{po}(x)}=1.$$

So we have proved that (a-1) and (a-2) are equivalent and imply
(a-3). To finish the proof, we will prove that (a-3) implies
(a-2). To that end it suffices with proving that
$$\lim_{x\to\infty}\frac{\overline{po}(x)-\overline{po}(x+y)}{\overline{po}(x)}=0,\qquad
\text{for any}\ y>0.$$ Indeed, using the equation (EA) we have
that for any $y>0,$
\begin{equation*}
\begin{split}
\mu-\frac{\int^{\infty}_0\mathrm{d}z \overline{\Pi}^+(z+x)}{\overline{po}(x)}&=\frac{\int^{\infty}_0\mathrm{d}z\overline{ne}(z)\overline{po}(x)+\widehat{d}\overline{po}(x)-\int^{\infty}_x po(\mathrm{d}y)\int^{y-x}_0\mathrm{d}z ne(z) -\widehat{d}\overline{po}(x)}{\overline{po}(x)}\\
&=\int^\infty_0\mathrm{d}z \overline{ne}(z)\frac{\overline{po}(x)-\overline{po}(x+z)}{\overline{po}(x)}\\
&\geq \int^\infty_y\mathrm{d}z \overline{ne}(z)\frac{\overline{po}(x)-\overline{po}(x+y)}{\overline{po}(x)}\geq 0.
\end{split}
\end{equation*}
and the assertion follows making $x\to \infty$ in the latter
equation since by assumption its left hand term tends to 0 as
$x\to\infty.$\end{proof}

\begin{proof}[Proof of (b) in Theorem~\ref{chnasc:1}]
The assumption that $\widehat{k}_0>0,$ implies that the renewal measure $\widehat{V}(\mathrm{d}y)$ is a finite measure and $\widehat{V}[0,\infty[=1/\widehat{k}_0.$ Thus if  $\overline{\Pi}^+(x)\in \mathcal{L}^0$ we have by the equation (EA) and the dominated convergence theorem that
\begin{equation*}
\begin{split}
\lim_{x\to\infty}\frac{\overline{po}(x)}{\overline{\Pi}^+(x)}&=\lim_{x\to\infty} \int^\infty_0\widehat{V}(\mathrm{d}y)\frac{\overline{\Pi}^+(x+y)}{\overline{\Pi}^+(x)}\\
&=\int^\infty_0\widehat{V}(\mathrm{d}y)\lim_{x\to\infty}\frac{\overline{\Pi}^+(x+y)}{\overline{\Pi}^+(x)}=\frac{1}{\widehat{k}_0}.
\end{split}
\end{equation*}

Now, that~(\ref{taub:s3}) holds is a straightforward consequence of the following identity, for any $x>0$
\begin{equation}\label{talacha}
\begin{split}
\overline{\Pi}^+(x)&=\int^{\infty}_0 ne(\mathrm{d}y)\left(\overline{po}(x)-\overline{po}(x+y)\right)+\widehat{k}_0 \overline{po}(x) + \widehat{d}\overline{p}(x)\\
&=\overline{po}(x)\int^1_{0} ne(\mathrm{d}y)\left(\frac{\overline{po}(x)-\overline{po}(x+y)}{\overline{po}(x)}\right)+ \int^{\infty}_1 ne(\mathrm{d}y)\left(\overline{po}(x)-\overline{po}(x+y)\right)\\
&\quad + \widehat{k}_0 \overline{po}(x) +\widehat{d}\overline{p}(x),
\end{split}
\end{equation}
which is obtained using the equation (EAI) and Fubini's theorem.
We have so proved that (b-1) implies (b-2) and (b-3). Next, to
prove that (b-2) implies (b-1) and (b-3) we assume that $po\in
\mathcal{L}^0$ and (\ref{taub:s3}) holds. Under this assumptions
we claim that $$\overline{\Pi}^+(x)\sim \widehat{k}_0
\overline{po}(x)+\widehat{d}\overline{p}(x) \qquad\text{as}\
x\to\infty.$$ Indeed, this can be deduced from
equation~(\ref{talacha}), using that
$\int^0_{-\infty}ne(\mathrm{d}y)\min\{\vert y \vert, 1\} <
\infty,$ that
$\lim_{x\to\infty}\overline{po}(x+y)/\overline{po}(x)=1$ for any
$y>0,$ and the dominated convergence theorem. Furthermore, we have
by hypothesis that $\widehat{d}\overline{p}(x)/\overline{po}(x)\to
0$ as $x\to\infty,$ which implies that
$$\overline{\Pi}^+(x)\sim \widehat{k}_0 \overline{po}(x)\qquad\text{as}\  x\to\infty.$$ To finish we next prove that (b-3) implies (b-2). Indeed, using the equation~(EAI) and the hypothesis (b-3) we get that
\begin{equation*}
\begin{split}
\frac{\overline{\Pi}^+(x)}{\overline{po}(x)}-\widehat{k}_0&=\frac{\int^{\infty}_x po(\mathrm{d}z)\overline{ne}(z-x)+\widehat{d}\overline{p}(x)}{\overline{po}(x)}\\
&=\int^{\infty}_0 ne(\mathrm{d}z)\frac{\left(\overline{po}(x)-\overline{po}(x+z)\right)}{\overline{po}(x)}+\frac{ \widehat{d}\overline{p}(x)}{\overline{po}(x)}\longrightarrow 0 \qquad\text{as} \ x\to \infty.
\end{split}
\end{equation*}
We deduce therefrom that (\ref{taub:s3}) holds and that $po\in\mathcal{L}^0$ since for any $y>0,$ $$\int^{\infty}_0 ne(\mathrm{d}z)\left(\frac{\overline{po}(x)-\overline{po}(x+z)}{\overline{po}(x)}\right)\geq \overline{ne}(y)\left(\frac{\overline{po}(x)-\overline{po}(x+y)}{\overline{po}(x)}\right)\geq 0.$$
\end{proof}

\begin{proof}[Proof of (a) in Proposition~\ref{bhvat0}]
According to Theorem~6.3.2 in \cite{thesevigon} under these assumptions the measure $po$
has infinite total mass if and only if $lim_{x\to 0+}\int^1_x\overline{\Pi}^+(z)\mathrm{d}z=\infty$ and in this case
$$\overline{po}(x)\sim \frac{1}{\widehat{d}}\int^1_x\overline{\Pi}^+(z)\mathrm{d}z, \qquad \text{as}\  x \to 0+.$$
Thus the assertion in (a) Proposition~\ref{bhvat0} is a consequence of this fact and the monotone density theorem for regularly varying functions. \end{proof}

\begin{proof}[Proof of (b) in Proposition~\ref{bhvat0}]
According to  Theorem~6.3.1 in \cite{thesevigon} under these assumptions, if we suppose $\lim_{x\to0+}\overline{po}(x)=\infty,$ then
$$\overline{po}(x)\sim \frac{1}{ne]0,\infty[+\widehat{k}_0}\overline{\Pi}^+(x),\qquad \text{as} \ x\to 0+.$$ The result follows.
\end{proof}

\end{subsection}

\begin{subsection}{Proof of Proposition~\ref{keythm}}
\begin{proof}[Sketch of proof of Proposition~\ref{keythm}]
The proof of the assertion (i) implies (ii) is a reworking of the
proof of Lemma~3 in Asmussen et al.~\cite{asmussenetal2002}, this
can be done in our setting since the only hypothesis needed in
that proof is that the dual ladder height has a finite mean.

To show that (i) implies (ii) in Proposition~\ref{keythm} we first
prove that under the assumption $\er(\widehat{H}_1)=\mu<\infty,$
the condition $\overline{\Pi}^+\in\mathcal{L}^0,$ implies that for
any $z>0,$
\begin{itemize}
\item[(BRT)]$$\displaystyle po]x,x+z[\sim
\frac{z}{\mu}\overline{\Pi}^+(x), \qquad
x\to\infty.$$\end{itemize} The latter estimate and the fact that
$\Pi^+\in\mathcal{L}^0$ implies that for any $a\geq 0,$
\begin{equation}\label{subexptails}\displaystyle po]x+a,x+a+z[\sim
\frac{z}{\mu}\overline{\Pi}^+(x), \qquad x\to\infty.\end{equation}
To prove that (BRT) holds, we may simply repeat the argument in
the proof of Lemma~3 in Asmussen et al.~\cite{asmussenetal2002}
using instead of the equation~(12) therein, the equation
$$po]x,x+z[=\int^\infty_x\Pi(dy)\widehat{V}]y-x-z,y-x[,\qquad
z>0,$$ which is an elementary consequence of equation (EAI) and
Fubini's theorem.

The result in (ii) in Proposition~\ref{keythm} follows from (BRT)
in the same way that the Key renewal theorem is obtained from
Blackwell's renewal theorem using the estimate in
(\ref{subexptails}) and the bounds
$$\frac{po]x,x+z[}{\overline{\Pi}^+(x)}\leq \widehat{V}(z), \qquad x>0, z>0,$$ which are
a simple consequence of the former equation and the fact that
$\widehat{V}$ is a renewal measure and so that for any $0<z<y,$
$\widehat{V}(y)-\widehat{V}(y-z)\leq \widehat{V}(z).$

To show that (ii) implies (i) we have to verify that for any $a>0$
$$\lim_{x\to\infty}\frac{\overline{\Pi}^+(x+a)}{\overline{\Pi}^+(x)}=1.$$
This is indeed true since using (ii) in Proposition~\ref{keythm}
it is straightforward that for any $z>0$ the assertion in (BRT)
holds and a further application of (ii) in
Proposition~\ref{keythm} to the function
$g_a(\cdot)=1_{\{]a,a+1[\}}(\cdot),$ $a>0,$ gives that for any
$a>0,$
$$\lim_{x\to\infty}\frac{po]x+a,x+a+1[}{\overline{\Pi}(x)}=\frac{1}{\mu},$$
and therefore, for any $a>0,$
$$\lim_{x\to\infty}\frac{\overline{\Pi}^+(x+a)}{\overline{\Pi}^+(x)}=\lim_{x\to\infty}\frac{\overline{\Pi}^+(x+a)}{po]x+a,x+a+1[}\frac{po]x+a,x+a+1[}{\overline{\Pi}^+(x)}=1.$$
\end{proof}
\end{subsection}

\end{section}

\begin{gracias}
It is my pleasure to thank Jean Bertoin for pointing out to me the work~\cite{MR85e:60093} and for insightful discussions and Ron Doney to bringing to my attention the papers~\cite{MR0423543} and~\cite{Grubel85}, suggesting to me an early version of Theorem~\ref{chnasc:1} and Proposition~\ref{bhvat0} and for a number of helpful comments. Part of this paper was written while I was visiting the Mathematical Sciences Institute of the Australian National University (ANU) in March 2005. I heartily thank Claudia Kl\"uppelberg and Ross Maller for arranging this visit, their hospitality and valuable discussions. Last but not least my thanks go to Alex Lindner for his hospitality during my visit to ANU and many stimulating discussions.
\end{gracias}

\bibliography{pathtrssmp}

\begin{thebibliography}{10}

\bibitem{asmussenetal2002}
S.~Asmussen, V.~Kalashnikov, D.~Konstantinides, C.~Kl{\"u}ppelberg, and
  G.~Tsitsiashvili.
\newblock A local limit theorem for random walk maxima with heavy tails.
\newblock {\em Statist. Probab. Lett.}, 56(4):399--404, 2002.

\bibitem{MR98e:60117}
J.~Bertoin.
\newblock {\em L\'evy processes}, volume 121 of {\em Cambridge Tracts in
  Mathematics}.
\newblock Cambridge University Press, Cambridge, 1996.

\bibitem{MR96b:60190}
J.~Bertoin and R.~A. Doney.
\newblock Cram\'er's estimate for {L}\'evy processes.
\newblock {\em Statist. Probab. Lett.}, 21(5):363--365, 1994.

\bibitem{MR90i:26003}
N.~H. Bingham, C.~M. Goldie, and J.~L. Teugels.
\newblock {\em Regular variation}, volume~27 of {\em Encyclopedia of
  Mathematics and its Applications}.
\newblock Cambridge University Press, Cambridge, 1989.

\bibitem{Deweert}
F.~De~{W}eert.
\newblock Attraction to stable distributions for l{\'e}vy processes at zero.
\newblock Technical report, University of {M}anchester, 2003.

\bibitem{doney2001}
R.~Doney.
\newblock Fluctuation theory for {L}\'evy processes.
\newblock In {\em L\'evy processes}, pages 57--66. Birkh\"auser Boston, Boston,
  MA, 2001.

\bibitem{MR2003d:60090}
R.~A. Doney and R.~A. Maller.
\newblock Stability of the overshoot for {L}\'evy processes.
\newblock {\em Ann. Probab.}, 30(1):188--212, 2002.

\bibitem{dynkin55}
E.~B. Dynkin.
\newblock Some limit theorems for sums of independent random quantities with
  infinite mathematical expectations.
\newblock {\em Izv. Akad. Nauk SSSR. Ser. Mat.}, 19:247--266, 1955.

\bibitem{Fristedt74}
B.~Fristedt.
\newblock Sample functions of stochastic processes with stationary, independent
  increments.
\newblock In {\em Advances in probability and related topics, Vol. 3}, pages
  241--396. Dekker, New York, 1974.

\bibitem{MR45:1250}
B.~E. Fristedt and W.~E. Pruitt.
\newblock Lower functions for increasing random walks and subordinators.
\newblock {\em Z. Wahrscheinlichkeitstheorie und Verw. Gebiete}, 18:167--182,
  1971.

\bibitem{MR89a:26002}
J.~L. Geluk and L.~de~Haan.
\newblock {\em Regular variation, extensions and {T}auberian theorems},
  volume~40 of {\em CWI Tract}.
\newblock Stichting Mathematisch Centrum Centrum voor Wiskunde en Informatica,
  Amsterdam, 1987.

\bibitem{gihmanskorohod}
{\u{I}}.~{\=I}. G{\={\i}}hman and A.~V. Skorohod.
\newblock {\em The theory of stochastic processes. {II}}.
\newblock Springer-Verlag, New York, 1975.
\newblock Translated from the Russian by Samuel Kotz, Die Grundlehren der
  Mathematischen Wissenschaften, Band 218.

\bibitem{MR85e:60093}
P.~Greenwood, E.~Omey, and J.~L. Teugels.
\newblock Harmonic renewal measures.
\newblock {\em Z. Wahrsch. Verw. Gebiete}, 59(3):391--409, 1982.

\bibitem{Grubel85}
R.~Gr{\"u}bel.
\newblock Tail behaviour of ladder-height distributions in random walks.
\newblock {\em J. Appl. Probab.}, 22(3):705--709, 1985.

\bibitem{kluppelberg88}
C.~Kl{\"u}ppelberg.
\newblock Subexponential distributions and integrated tails.
\newblock {\em J. Appl. Probab.}, 25(1):132--141, 1988.

\bibitem{KKM}
C.~Kl{\"u}ppelberg, A.~E. Kyprianou, and R.~A. Maller.
\newblock Ruin probabilities and overshoots for general {L}\'evy insurance risk
  processes.
\newblock {\em Ann. Appl. Probab.}, 14(4):1766--1801, 2004.

\bibitem{rogozin71}
B.~A. Rogozin.
\newblock Distribution of the first laddar moment and height, and fluctuations
  of a random walk.
\newblock {\em Teor. Verojatnost. i Primenen.}, 16:539--613, 1971.

\bibitem{sato99}
K.-I. Sato.
\newblock {\em L\'evy processes and infinitely divisible distributions},
  volume~68 of {\em Cambridge Studies in Advanced Mathematics}.
\newblock Cambridge University Press, Cambridge, 1999.

\bibitem{sinai57}
Y.~G. Sina\v\i.
\newblock On the distribution of the first positive sum for a sequence of
  independent random variables.
\newblock {\em Theory Prob.}, 16:122--129, 1957.

\bibitem{MR0423543}
N.~Veraverbeke.
\newblock Asymptotic behaviour of {W}iener-{H}opf factors of a random walk.
\newblock {\em Stochastic Processes Appl.}, 5(1):27--37, 1977.

\bibitem{thesevigon}
V.~Vigon.
\newblock {\em Simplifiez vos {L\'e}vy en titillant la factorisation de
  {W}iener-{H}opf}.
\newblock PhD thesis, Universit\'e {L}ouis {P}asteur, 2002.

\bibitem{MR2002i:60101}
V.~Vigon.
\newblock Votre {L}\'evy rampe-t-il?
\newblock {\em J. London Math. Soc. (2)}, 65(1):243--256, 2002.

\end{thebibliography}
\end{document}